\magnification=\magstep1
%%%%%%%%%%%%%%%%%%%%%%%%%%%%%%%%%%%%%%%%%%%%%%%%%%%%%%%%%%%%%%%%%%%%
%%
%%  PAPERMAC.TEX
%%
%%  (c) 1986 D.E. Taylor
%%
%%  20 November 85     23 February 87
%%
%%%%%%%%%%%%%%%%%%%%%%%%%%%%%%%%%%%%%%%%%%%%%%%%%%%%%%%%%%%%%%%%%%%%
\magnification=\magstep1
%
%  font section
%
\font\titlefont=cmbx10 scaled \magstep1
\font\secfont  =cmbx10
               % sans serif
\font    \affil=cmsl8                % small slanted
\font \headfont=cmcsc10              % caps and small caps
\font \footfont=cmr8
\font      \csc=cmcsc8               % caps and small caps

%  page shape

\baselineskip=12pt
\hoffset=1truecm
%\voffset=1truecm
\hsize=11.667cm
\vsize=18.333cm
\parskip=5pt plus 2pt
\parindent=0pt

%  running titles and other headings
%
\newif\iftitle
\newcount\tagno   % for use by the calling file to number equations etc.
\def\papername{}
\def\title#1\par{\global\titletrue{\baselineskip=24pt
   \halign{\line{\titlefont\hfil##\hfil}\cr#1\unskip\cr}}
   \message{#1}\tagno=0\vskip36pt\noindent
   \centerline{D. E. Taylor}
   \medskip
   \centerline{\affil Department of Pure Mathematics}
   \centerline{\affil The University of Sydney}
   \centerline{\affil Australia 2006}
   \vskip 42pt\noindent
   \xdef\papername{#1\unskip}}

\def\shorttitle#1{\xdef\papername{#1\unskip}\message{#1}\noindent}

\headline={\iftitle\hfil
   \else\hfil\headfont\papername\hfil\fi}
\footline={\iftitle\global\titlefalse\hfill
   \else\hfil\rm\folio\hfil\fi}
\outer\def\entry#1\par{\medbreak\vskip -\parskip\noindent
   \global\advance\tagno by 1\hbox{\bf\the\tagno.\quad}{\sl #1}
   \nobreak\medskip\noindent\message{(\the\tagno)}}
%
%   SYMBOLS
%
\def\GF{{\mit GF}}
\def\GL{{\mit GL}}
\def\SL{{\mit SL}}
\def\Sp{{\mit Sp}}
\def\SU{{\mit SU}}

\def\bdots{\mathinner{\mskip1mu\vbox{\kern7pt\hbox{.}}\mskip2mu
   \raise3pt\hbox{.}\mskip2mu\raise6pt\hbox{.}\mskip1mu}}
%
%  \fbox from LaTeX
%
\chardef\prevatcode=\catcode`\@       \catcode`\@=11
\newdimen\fboxrule    \fboxrule= .4\p@
\newdimen\fboxsep     \fboxsep= 4\p@
\def\fbox#1{\leavevmode	
   \setbox0\hbox{#1}\dimen@=\fboxrule
   \advance\dimen@ \fboxsep
   \advance\dimen@ \dp0
   \hbox{\lower\dimen@\hbox
      {\vbox{\hrule height \fboxrule
         \hbox{\vrule width \fboxrule\hskip\fboxsep
               \vbox{\vskip\fboxsep\box0\vskip\fboxsep}%
         \hskip\fboxsep\vrule width \fboxrule}%
         \hrule height \fboxrule}}}}
\catcode`\@=\prevatcode
%
%   REFERENCES
%
\def\references{\medbreak\bigskip
   \leftline{\secfont References}\nobreak\smallskip\noindent
   \begingroup
   \parindent=0pt
   \pretolerance=300\tolerance=500\finalhyphendemerits=0
   \parskip=3pt plus 1pt minus 1pt\parindent=12pt
   \spaceskip=4pt plus 2pt minus 1pt
   \xspaceskip=5pt plus 3pt minus 2pt 
   \sfcode `\.=1000\sfcode `\,=2000}
\def\endref{\endgroup}

\def\medfil{\par\vfil\penalty-50\vfilneg}
\title Pairs of Generators for Matrix Groups. I

\shorttitle {Generators for Matrix Groups}

It has been shown by Steinberg (1962) that every finite simple 
group of Lie type can be generated by two elements.  These 
groups can be constructed from the simple Lie algebras over the
complex numbers by the methods of Chevalley (1955), Steinberg (1959)
and Ree (1961).  The generators obtained by Steinberg (1962) are 
given in terms of the root structure of the corresponding Lie 
algebras.

The identification of the groups of Lie type $A_n$, $B_n$, $C_n$ and
$D_n$ with classical matrix groups is due to Ree (1957).
and an exposition of his results can be found in the book 
of Carter (1972).  The proofs ultimately rely on the
work of Dickson (1901).  In this note we give tables of 
generators for the groups $\GL(n,q)$, $\SL(n,q)$, $\Sp(2n,q)$,
$U(n,q)$ and $\SU(n,q)$.  For the most part, the generators have 
been obtained by translating Steinberg's generators into matrix 
form via the methods of Ree (1957)\footnote{}{\footfont The Cayley Bulletin
No.~3, October 1987, 76--85.}.
\medskip

{\bf Notation}

Let $E_{ij}$ denote a square matrix with 1 in the $(i,j)$th
position and 0 elsewhere.  For $\alpha\in\GF(q)$ and 
$i\ne j$ we set
$$x_{ij}(\alpha) = I + \alpha E_{ij}$$
and we let $h_i(\alpha)$ denote the diagonal matrix obtained by 
replacing the $i$th entry of the identity matrix by $\alpha$.
The $x_{ij}(\alpha)$ are the root elements of $\SL(n,q)$.

Let $w_i$ denote the monomial matrix obtained from the 
permutation matrix corresponding to the transposition $(i,i+1)$ 
by replacing the $(i+1,i)$-th entry by $-1$.  Then
$w=w_1w_2\ldots w_{n-1}$ represents the $n$-cycle $(1,2,\ldots,n)$.

Let $\xi$ be a generator of the multiplicative 
group of $\GF(q)$.  
\medfil
\entry $\GL(n,q)$, $q\ne 2$

Generators for $\GL(n,q)$ are
$$h_1(\xi) = \pmatrix{
    \xi & 0 & \ldots & 0\cr
      0 & 1 & \ldots & 0\cr
        &   & \ddots &  \cr
      0 & 0 & \ldots & 1\cr}\quad\hbox{and}\quad
x_{12}(1)w = \pmatrix{
 -1 & 0 & \ldots & 0 & 1\cr
 -1 & 0 & \ldots & 0 & 0\cr
  0 & -1& \ldots & 0 & 0\cr
    &   & \ddots &   &  \cr
  0 & 0 & \ldots & -1& 0\cr}.$$
When $n=2$, the generators are
$$\pmatrix{\xi&0\cr0&1\cr}\quad\hbox{and}\quad
\pmatrix{-1&1\cr-1&0\cr}.$$
Generators for $\GL(n,2) = \SL(n,2)$ are given below.
\bigskip

\entry $\SL(n,q)$, $q>3$

The generators are
$$h_1(\xi)h_2(\xi^{-1}) = \pmatrix{
\xi&&&&\cr&\xi^{-1}&&&\cr&&1&&\cr&&&\ddots&\cr&&&&1\cr}$$
and the matrix $x_{12}(1)w$ given above.
\medfil

\entry $\SL(n,2)$ and $\SL(n,3)$

The generators are
$$x_{12}(1) = \pmatrix{
1&1&0&\ldots&0\cr
0&1&0&\ldots&0\cr
0&0&1&\ldots&0\cr
&&&\ddots&\cr
0&0&0&\ldots&1\cr}\quad\hbox{and}\quad
w=\pmatrix{
0 & 0 & \ldots & 0 & 1\cr
-1& 0 & \ldots & 0 & 0\cr
0 &-1 & \ldots & 0 & 0\cr
  &   & \ddots &   &  \cr
0 & 0 & \ldots &-1 & 0\cr}.$$
\medfil
\entry $\Sp(2n,q)$, $q$ odd, $n>1$

The symplectic group  $\Sp(2n,q)$ consists of the $2n\times 2n$ 
matrices $X$ which satisfy the condition 
$$X^tJX = J,$$
where
$$J=\pmatrix{&&&\smash{\vrule height1ex depth72pt}&&&1\cr
             &&&&&\bdots&\cr
             &&&&1&&\cr
      \multispan7\hrulefill\cr
             &&\llap{$-$}1&&&&\cr
             &\bdots&&&&&\cr
             -1&&&&&&\cr}.$$

For $1\le i\le n$, let $i' = 2n-i+1$.  Define
$$\hat h_i(\alpha)=h_i(\alpha)h_{i'}(\alpha^{-1})$$
and 
$$\hat x_{ij}(\alpha) = x_{ij}(\alpha)x_{j'i'}(-\alpha).$$
Let $\hat w$ be the monomial matrix obtained from the permutation
matrix of the $2n$-cycle
$$(1,2,\ldots,n,1',2',\ldots,n')$$
by replacing the $(2n,n)$th entry by $-1$.

\medfil
Generators for $\Sp(2n,q)$, $q$ odd, are
$$\hat h_1(\xi)= \pmatrix{
\xi &&&&\smash{\vrule height1ex depth108pt}&&&&\cr
       &1&&&   &&&&\cr
    &&\ddots&& &&&&\cr
       &&&1&   &&&&\cr
      \multispan9\hrulefill\cr
       &&&&    &1&&&\cr
       &&&&    &&\ddots&&\cr
       &&&&    &&&1&\cr
       &&&&    &&&&\xi^{-1}\cr}$$
\medfil
and
$$\hat x_{12}(1)\hat w = \pmatrix{
 1&0&&&&\smash{\vrule height1ex depth132pt}&1&&&&\cr
 1&0&&&&       &&&&&\cr
 0&1&&&&       &&&&&\cr
 &&\ddots&&&   &&&&&\cr
 &&&1&0&       &0&&&&\cr
\multispan{11}\hrulefill\cr
 &&&0&0&       &0&1&&&\cr
 &&&&&         &&0&1&&\cr
 &&&&&         &&&\smash{\ddots}&&\cr
 &&&0&1&       &&&&0&1\cr
 &&&0&\llap{$-$}1&       &&&&&0\cr}.$$
\medfil
When $n=2$, the matrices are
$$\pmatrix{\xi&0&0&0\cr
                0&1&0&0\cr
                0&0&1&0\cr
                0&0&0&\xi^{-1}\cr}\quad
\hbox{and}\quad\pmatrix{
1&0&1&0\cr
1&0&0&0\cr
0&1&0&1\cr
0&\llap{$-$}1&0&0\cr}.$$

Note that $\Sp(2,q)\simeq\SL(2,q)$.

\medfil
\entry $\Sp(2n,q)$, $q$ even, $q\ne 2$, $n>1$

The $\hat x_{ij}(\alpha)$ are the short root elements of $\Sp(2n,q)$.
The long root elements are transvections
$\hat z_i(\alpha)=x_{ii'}(\alpha)$.  

Generators for $\Sp(2n,q)$ are
$$\hat h_1(\xi)\hat h_n(\xi) = \pmatrix{
\xi &&&&\smash{\vrule height1ex depth108pt}&&&&\cr
       &1&&&        &&&&\cr
    &&\ddots&&      &&&&\cr
       &&&\xi&   &&&&\cr
      \multispan9\hrulefill\cr
       &&&&         &\xi^{-1}&&&\cr
       &&&&         &&\ddots&&\cr
       &&&&         &&&1&\cr
       &&&&         &&&&\xi^{-1}\cr}$$
\medfil
and
$$\hat x_{1n}(1)\hat z_1(1)\hat w =\pmatrix{
 0&0&&1&1&\smash{\vrule height1ex depth132pt}&1&&&&\cr
 1&0&&0&0&        &0&&&&\cr
 0&1&&&&          &&&&&\cr
 &&\ddots&&&      &&&&&\cr
 &&&1&0&          &0&&&&\cr
\multispan{11}\hrulefill\cr
 &&&&1&         &0&1&&&\cr
 &&&&0&         &&0&1&&\cr
 &&&&&          &&&\smash{\ddots}&&\cr
 &&&&&          &&&&0&1\cr
 &&&&1&         &&&&&0\cr}.$$
\medfil
For $\Sp(4,q)$ the matrices become
$$\pmatrix{
\xi&0&0&0\cr
0&\xi&0&0\cr
0&0&\xi^{-1}&0\cr
0&0&0&\xi^{-1}\cr}\quad\hbox{and}\quad
\pmatrix{
1&1&1&0\cr
1&0&0&0\cr
0&1&0&1\cr
0&1&0&0\cr}.$$
\medfil
\entry $\Sp(2n,2)$, $n>2$

The generators are
$$\hat x_{1n}(1)\hat z_1(1) = \pmatrix{
1&&1&\smash{\vrule height1ex depth84pt}&&&1\cr
&\ddots&&   &&&\cr
&&1&        &&&\cr
\multispan7\hrulefill\cr
&&&         &1&&1\cr
&&&         &&\ddots&\cr
&&&         &&&1\cr}$$
\medfil
and
$$\hat w = \pmatrix{
 0&&&&&\smash{\vrule height1ex depth132pt}&1&&&&\cr
 1&0&&&&        &&&&&\cr
 &&\ddots&&&    &&&&&\cr
 &&&0&&         &&&&&\cr
 &&&1&0&        &0&&&&\cr
\multispan{11}\hrulefill\cr
 &&&&0&         &0&1&&&\cr
 &&&&&          &&0&&&\cr
 &&&&&          &&&\ddots&&\cr
 &&&&0&         &&&&0&1\cr
 &&&&1&         &&&&&0\cr}.$$
\medfil
\entry $\Sp(4,2)$

The matrices are
$$\pmatrix{
1&0&1&1\cr
1&0&0&1\cr
0&1&0&1\cr
1&1&1&1\cr}\quad\hbox{and}\quad
\pmatrix{
0&0&1&0\cr
1&0&0&0\cr
0&0&0&1\cr
0&1&0&0\cr}$$
\medfil
\entry $U(2n,q)$, $n>1$

If $x\in \GF(q^2)$, we set $\bar x=x^q$.  If $X$ is a matrix, $\overline X$ is
obtained from $X$ by replacing each entry $x$ with $\bar x$.  The
unitary group consists of the matrices $X$ such that
$$\overline X^tJX = J,$$
where
$$J=\pmatrix{&&1\cr
             &\bdots&\cr
             1&&\cr}.$$

Let $\xi$ be a primitive element of $\GF(q^2)$ and let $\eta$
be an element of trace 0, i.e., $\eta + \overline\eta = 0$.
If $q$ is odd, we may take 
$\eta = \xi^{(q+1)/2}$.  If $q$ is even, we may take $\eta = 1$.
\medfil
For this section define $i'= 2n+1-i$ and set
$$\tilde h_i(\alpha)= h_i(\alpha)h_{i'}(\overline\alpha^{-1}),$$
and
$$\tilde x_{ij}(\alpha)=x_{ij}(\alpha)x_{j'i'}(-\overline\alpha),$$
for $1\le i,j\le n$.
\medfil
The matrix $\tilde w$ is similar to $\hat w$ of previous 
sections except that here it corresponds to the permutation
$(1,2,\ldots,n,1',2',\ldots,n')$ and it has $\eta$ in the
$(1,n+1)$th position and $-\eta^{-1}$ in the $(2n,n)$th 
position.
\medfil
Generators for $U(2n,q)$ are
$$\tilde h_1(\xi)=\pmatrix{
\xi &&&&\smash{\vrule height1ex depth108pt}&&&&\cr
       &1&&&   &&&&\cr
    &&\ddots&& &&&&\cr
       &&&1&   &&&&\cr
      \multispan9\hrulefill\cr
       &&&&    &1&&&\cr
       &&&&    &&\ddots&&\cr
       &&&&    &&&1&\cr
       &&&&    &&&&\bar\xi^{-1}\cr}$$
\medfil
and
$$\tilde x_{12}(1)\tilde w =\pmatrix{
 1&0&&&&\smash{\vrule height1ex depth132pt}&\eta&&&&\cr
 1&0&&&&       &&&&&\cr
 0&1&&&&       &&&&&\cr
 &&\ddots&&&   &&&&&\cr
 &&&1&0&       &0&&&&\cr
\multispan{11}\hrulefill\cr
 &&&&0&        &0&1&0&&\cr
 &&&&&         &&0&1&&\cr
 &&&&&         &&&\ddots&&\cr
 &&&&\eta^{-1}&       &&&&0&1\cr
 &&&&\llap{$-$}\eta^{-1}&       &&&&&0\cr}.$$
\medfil
When $n=2$ the matrices are
$$\pmatrix{
\xi&0&0&0\cr
   0&1&0&0\cr
   0&0&1&0\cr
   0&0&0&\bar\xi^{-1}\cr}
\quad\hbox{and}\quad
\pmatrix{
1&0&\eta&0\cr
1&0&0&0\cr
0&\eta^{-1}&0&1\cr
0&\llap{$-$}\eta^{-1}&0&0\cr}.$$
\medfil
\entry $\SU(2n,q)$, $n>1$

Generators are
$$\tilde h_1(\xi)\tilde h_2(\xi^{-1})=\pmatrix{
\xi &&&&\smash{\vrule height1ex depth108pt}&&&&\cr
       &\xi^{-1}&&&   &&&&\cr
    &&\ddots&& &&&&\cr
       &&&1&   &&&&\cr
      \multispan9\hrulefill\cr
       &&&&    &1&&&\cr
       &&&&    &&\ddots&&\cr
       &&&&    &&&\bar\xi&\cr
       &&&&    &&&&\bar\xi^{-1}\cr}$$
and $\tilde x_{12}(1)\tilde w$.
\medfil
\entry $U(2n+1,q)$

In this section let $i' = 2n+2-i$ and define $\tilde h_i(\alpha)$
and $\tilde x_{ij}(\alpha)$ as before.  In the following matrices the
boxed entry is in position $(n+1,n+1)$.
\medfil
For $\alpha,\beta\in \GF(q^2)$ such that 
$\alpha\overline\alpha + \beta +\overline\beta = 0$, set
$$Q(\alpha,\beta)=\pmatrix{
I&&&&\cr
 &1&\alpha&\beta&\cr
 &&\fbox{1}&\llap{$-$}\overline\alpha&\cr
 &&&1&\cr
 &&&&I\cr}.$$

Let $w'$ be the monomial matrix obtained from the permutation matrix
of $(n',\ldots,2',1',n,\ldots,2,1)$ by replacing the $(n+1,n+1)$st entry
by $-1$.
\medfil
Let $\beta$ be an element of $\GF(q^2)$ such that $\beta+\overline \beta= -1$.
We may take $\beta = -(1+\bar\xi/\xi)^{-1}$.
\medfil
Generators are
$$\tilde h_n(\xi)=\pmatrix{
      1&&&&&&\cr
       &\ddots&&&&&\cr
       &&\xi&&&&\cr
       &&&\fbox{1}&&&\cr
       &&&&\bar\xi^{-1}&&\cr
       &&&&&\ddots&\cr
       &&&&&&1\cr}$$
\medfil
and 
$$Q(1,\beta)w'=\pmatrix{
      0&1&&&&&&&&&\cr
       &0&&&&&&&&&\cr
       &&\ddots&&&&&&&&\cr
       &&&0&1&&&&&&\cr
     \beta &&&&0&-1&&&&&1\cr
     -1&&&&&\fbox{$-1$}&&&&&\cr
      1&&&&&&0&&&&\cr
       &&&&&&1&0&&&\cr
       &&&&&&&&\ddots&&\cr
       &&&&&&&&&0&\cr
       &&&&&&&&&1&0\cr}$$
\medfil
\entry $\SU(2n+1,q)$, $n\ne 1$ or $q\ne 2$

Generators are 
$$\tilde h_n(\xi)\tilde h_{n+1}(\xi^{-1})=\pmatrix{
      1&&&&&&\cr
       &\ddots&&&&&\cr
       &&\xi&&&&\cr
       &&&\fbox{$\bar\xi/\xi$}&&&\cr
       &&&&\bar\xi^{-1}&&\cr
       &&&&&\ddots&\cr
       &&&&&&1\cr}$$
and $Q(1,\beta)w'$ as above.
\medfil
\entry $\SU(3,2)$

Generators are
$$\pmatrix{
1&\xi&\xi\cr
0&1&\xi^2\cr
0&0&1\cr}\quad\hbox{and}\quad
\pmatrix{
\xi&1&1\cr
1&1&0\cr
1&0&0\cr}.$$

\medskip
{\bf Cayley Functions}

The groups discussed in this note are defined in Cayley\footnote*{\footfont
In 1993 Cayley was replaced by {\csc Magma}} with these generators
by means of the following functions:

\bigskip
\centerline{$\GL(n,q)$\enspace -- -- -- {\it general linear$(n,q)$}}

\bigskip
\centerline{$\SL(n,q)$\enspace -- -- -- {\it special linear$(n,q)$}}

\bigskip
\centerline{$\Sp(2n,q)$\enspace -- -- -- {\it symplectic$(2 * n,q)$}}

\bigskip
\centerline{$U(n,q)$\enspace -- -- -- {\it general unitary$(n,q)$}}

\bigskip
\centerline{$\SU(n,q)$\enspace -- -- -- {\it special unitary$(n,q)$}}

\references
R. W. Carter, {\sl `Simple Groups of Lie Type,'} 
(Wiley-Interscience: New-York 1972).

C. Chevalley, Sur certaines groupes simples, {\sl T\^ohoku Math.
J. \bf 7} (1955), 14-66.

L. E. Dickson, {\sl `Linear groups,'} (Teubner: Leipzig 1901).

Rimhak Ree, On some simple groups defined by C. Chevalley,
{\sl Trans. Amer. Math. Soc. \bf 84} (1957), 392-400.

Rimhak Ree, A family of simple groups associated with the simple
Lie algebras of type $(F_4)$, {\sl Amer. J. Math. \bf 83} (1961), 401-420.

Rimhak Ree, A family of simple groups associated with the simple
Lie algebras of type $(G_2)$, {\sl Amer. J. Math. \bf 83} (1961), 432-462.

Robert Steinberg, Automorphisms of finite linear groups,
{\sl Canad. J. Math. \bf 12} (1960), 606-615.

Robert Steinberg, Generators for simple groups, {\sl
Canad. J. Math. \bf 14} (1962), 277-283.
\endref
\bye